\documentclass[10pt]{article}

\usepackage{graphicx}%
\usepackage{amsmath,amssymb,amsthm,upref,bm}%

\newtheorem{theorem}{Theorem}%
\begin{document}

\baselineskip=4.4mm

\makeatletter

\newcommand{\E}{\mathrm{e}\kern0.2pt} 
\newcommand{\D}{\mathrm{d}\kern0.2pt}
\newcommand{\RR}{\mathbb{R}}
\newcommand{\CC}{\mathbb{C}}%
\newcommand{\ii}{\kern0.05em\mathrm{i}\kern0.05em}

\renewcommand{\Re}{\mathrm{Re}} 
\renewcommand{\Im}{\mathrm{Im}}

\def\bottomfraction{0.9}

\title{\bf Inverse mean value property of solutions \\ to the modified Helmholtz
equation}

\author{Nikolay Kuznetsov}

\date{}

\maketitle

\vspace{-8mm}

\begin{center}
Laboratory for Mathematical Modelling of Wave Phenomena, \\ Institute for Problems
in Mechanical Engineering, Russian Academy of Sciences, \\ V.O., Bol'shoy pr. 61, St
Petersburg 199178, Russian Federation \\ E-mail: nikolay.g.kuznetsov@gmail.com
\end{center}

\begin{abstract}
\noindent A theorem characterizing analytically balls in the Euclidean space $\RR^m$
is proved. For this purpose positive solutions of the modified Helmholtz equation
are used instead of harmonic functions applied in previous results. The obtained
Kuran type theorem is based on the volume mean value property of solutions to this
equation.
\end{abstract}

\setcounter{equation}{0}

\section{Introduction and main result}

In 1972 Kuran \cite{K} proved the following inverse of the volume mean
value theorem for harmonic functions:
\begin{quote}
{\it Let $D$ be a domain (= connected open set) of finite (Lebesgue) measure in the
Euclidean space $\RR^m$ where $m \geq 2$. Suppose that there exists a point $P_0$
in $D$ such that, for every function $h$ harmonic in $D$ and integrable over $D$,
the volume mean of $h$ over $D$ equals $h (P_0)$. Then $D$ is an open ball (disk
when $m=2$) centred at $P_0$.}
\end{quote}
The result was originally obtained by Epstein \cite{E} for a simply connected
two-dimensional~$D$. Armitage and Goldstein \cite{AG} proved this result assuming
that the mean value equality holds only for positive harmonic functions which are
$L^p$-integrable, $p \in  (0, n / (n-2))$. Hansen and Netuka \cite{HN1}
considered some particular class of potentials as the set of test harmonic functions
in Kuran's theorem. A slight modification of his considerations shows that Kuran's
theorem is valid even if $D$ is disconnected; see \cite{NV}, p.~377.

In the survey article \cite{NV}, one finds also a discussion of applications of
Kuran's theorem and a possibility of similar results involving some kinds of average
over $\partial D$, where $D$ is a bounded domain. One of them (due to
Kosmodem'yanskii \cite{Kos}) is based on the relation similar to that between the
mean values over balls and spheres and reads as follows:
\begin{quote} 
{\it Let $D \subset \RR^2$ be a bounded, convex $C^2$-domain. If the equality
\[ \frac{1}{|D|} \int_{D} u (x) \, \D x = \frac{1}{|\partial D|} \int_{\partial D} 
u (x) \, \D S_x
\]
holds for every function $u \in C^2 (D) \cap C^1 (\overline D)$ which is harmonic
in $D$, then $D$ is an open disc.}
\end{quote}
Here and below $|D|$ is the domain's area (volume if $D \subset \RR^m$, $m \geq 3$),
whereas $|\partial D|$ is the boundary's length (area if $D \subset \RR^m$, $m \geq
3$), and $|B_r| = \omega_m r^m$ is the volume of a ball $B_r$ of radius $r$; the
volume of unit ball is $\omega_m = 2 \, \pi^{m/2} / [m \Gamma (m/2)]$, whereas
$\Gamma$ denotes the Gamma function.

In this note, we prove a new analytic characterization of balls. Like Kuran's
theorem, it is based on the $m$-dimensional volume mean value equality, but uses
solutions of the modified Helmholtz equation
\begin{equation}
\nabla^2 u - \lambda^2 u = 0 , \quad \lambda \in \RR \setminus \{0\} 
\label{Hh}
\end{equation}
instead of harmonic functions; here $\nabla = (\partial_1, \dots , \partial_m)$
denotes the gradient operator and $\partial_i = \partial / \partial x_i$. Solutions
are assumed to be real; indeed, the obtained results can be extended to
complex-valued functions by considering the real and imaginary part separately.

Before giving the precise formulation of the main result, let us introduce some
notation. By $B_r (x) = \{ y : |y-x| < r \}$ we denote the open ball of radius $r$
centred at $x \in \RR^m$; if $D \subset \RR^m$ is a bounded domain, then $D_r = D
\cup \left[ \cup_{x \in \partial D} B_r (x) \right]$ is its dilated copy such that
the distance from $\partial D_r$ to $D$ is equal to $r$. For a function $f$ integrable
over $D$, which has finite Lebesgue measure,
\[ M (f, D) = \frac{1}{|D|} \int_{D} f (x) \, \D x
\]
denotes its volume mean value over $D$. Also, we need the following function
\begin{equation}
a (t) = \Gamma \left(\frac{m}{2} + 1\right) \frac{I_{m/2} (t)}{(t / 2)^{m/2}} \, ,
\label{a}
\end{equation}
where $I_\nu$ stands for the modified Bessel function of order $\nu$. The relation
\begin{equation}
[z^{-\nu} I_\nu (z)]' = z^{-\nu} I_{\nu+1} (z) \ \ \mbox{(see \cite{Wa}, p. 79)} ,
\label{mon}
\end{equation}
where the right-hand side is positive for $z > 0$ and vanishes at $z = 0$, implies
that the function $a$ increases monotonically on $[0, \infty)$ from $a (0) = 1$ to
infinity; the latter is a consequence of the asymptotic formula valid as $|z| \to
\infty$:
\begin{equation*}
I_\nu (z) = \frac{\E^z}{\sqrt{2 \pi z}} \left[ 1 + O (|z|^{-1}) \right] \, , \ \
|\arg z| < \pi /2  \ \ \mbox{(see \cite{Wa}, p. 80)} \, .
\end{equation*}

The function $a$ arises in the $m$-dimensional mean value formula for balls
\begin{equation}
a (\lambda r) \, u (x) = \frac{1}{|B_r|} \int_{B_r (x)} u (y) \, \D y \, , \quad x
\in D , \label{MM}
\end{equation}
which holds, for example, if $u \in C^0 (\overline D)$ is a solution of \eqref{Hh}
in $D$ and $B_r (x) \subset D$. This equality was obtained by the author recently;
see \cite{Ku1}, p.~95. Before that only the three-dimensional mean value formula for
spheres had been derived by C.~Neumann (see his book \cite{NC}, Chapter~9,
Section~3, published in 1896), whereas the $m$-dimensional formula for spheres was
given without proof in \cite{Po}; its derivation see in the author's note
\cite{Ku2}.

Now, we are in a position to formulate the main result.

\begin{theorem}
Let $D \subset \RR^m$, $m \geq 2$, be a bounded domain such that its complement is
connected, and let $r$ be a positive number such that $|B_r| \leq |D|$. Suppose that
there exists a point $x_0 \in D$ such that for some $\lambda > 0$ the mean value
equality $u (x_0) \, a (\lambda r) = M (u, D)$ holds for every positive function $u$
satisfying equation \eqref{Hh} in $D_r$. If also $|D| = |B_r|$ provided $B_r (x_0)
\setminus \overline D \neq \emptyset$, then $D = B_r (x_0)$.
\end{theorem}

\section{Proof of Theorem 1 and discussion}

Prior to proving Theorem 1, we introduce the following function
\begin{equation}
U (x) = \Gamma \left(\frac{m}{2}\right) \frac{I_{(m-2)/2} (\lambda |x|)} {(\lambda
|x| / 2)^{(m-2)/2}} \, , \quad x \in \RR^m , \label{U}
\end{equation}
where the coefficient is chosen so that $U (0) = 1$. Let us consider some of its
properties. According to \eqref{mon}, this spherically symmetric function
monotonically increases as $|x|$ goes from zero to infinity. Also, it solves
equation \eqref{Hh} in $\RR^m$; indeed, the representation
\begin{equation}
U (x) = \frac{2 \, \Gamma (m/2)}{\sqrt \pi \, \Gamma ((m-1)/2)} \int_0^1 (1 -
s^2)^{(m-3)/2} \cosh (\lambda |x| s) \, \D s \, , \label{PU}
\end{equation}
is easy to differentiate, thus verifying \eqref{Hh}. This formula for $U$ is a
consequence of Poisson's integral (see \cite{NU}, p. 223):
\[ I_\nu (z) = \frac{(z / 2)^\nu}{\sqrt \pi \, \Gamma (\nu + 1/2)} \int_{-1}^1 (1 -
s^2)^{\nu - 1/2} \cosh z s \, \D s \, .
\]
Moreover, \eqref{PU} takes particularly simple form for $m=3$, namely, $U (x) =
(\lambda |x|)^{-1} \sinh \lambda |x|$. Since formulae \eqref{a} and \eqref{U} are
similar, Poisson's integral allows us to compare these functions. Indeed, the
inequality
\begin{equation} 
[ U (x) ]_{|x| = r} > a (\lambda r) \label{aU}
\end{equation}
immediately follows because $U (x)$ is spherically symmetric.

\begin{proof}[Proof of Theorem 1.]
Without loss of generality, we suppose that the domain $D$ is located so that $x_0$
coincides with the origin. Let us show that the assumption that $D \neq B_r (0)$
leads to a contradiction.

It is clear that either $B_r (0) \subset D$ or $B_r (0) \setminus \overline D \neq
\emptyset$ (the equality $|B_r| = |D|$ is assumed in the latter case), and we treat
these two cases separately. Let us consider the second case first and introduce the
bounded open sets $G_i = D \setminus \overline{B_r (0)}$ and $G_e = B_r (0)
\setminus \overline D$, for which we have $|G_e| = |G_i| \neq 0$ in view of
assumptions about $D$ and $r$. The volume mean equality for $U$ over $D$ can be
written as follows:
\begin{equation}
|D| \, a (\lambda r) = \int_D U (y) \, \D y \, ; \label{1}
\end{equation}
here the condition $U (0) = 1$ is taken into account. Since property \eqref{MM}
holds for $U$ over $B_r (0)$, we write it in the same way:
\begin{equation}
|B_r| \, a (\lambda r) = \int_{B_r (0)} U (y) \, \D y \, . \label{2}
\end{equation}
Subtracting \eqref{2} from \eqref{1}, we obtain
\begin{equation*}
0 = \int_{G_i} U (y) \, \D y - \int_{G_e} U (y) \, \D y > 0 \, .
\end{equation*}
Indeed, the difference is positive since $U (y)$ (positive and monotonically
increasing with~$|y|$) is greater than $[U (y)]_{|y| = r}$ in $G_i$ and less than
$[U (y)]_{|y| = r}$ in $G_e$, whereas $|G_i| = |G_e|$. This contradiction proves the
result in this case.

In the case when $B_r (0) \subset D$, we also have to obtain a contradiction when
$B_r (0) \neq D$. Now, after subtraction of \eqref{2} from \eqref{1} we have
\begin{equation*}
( |D| - |B_r| ) \, a (\lambda r) = \int_{G_i} U (y) \, \D y > |G_i| \, [ U (y)
]_{|y| = r} \, , \label{3}
\end{equation*}
where the last inequality is again a consequence of positivity of $U (y)$ and the
fact that it increases monotonically with $|y|$. It is clear that $|G_i| = |D| -
|B_r| > 0$ because $B_r (0) \subset~D$. Therefore, $a (\lambda r) > [ U (y) ]_{|y| =
r}$, which contradicts \eqref{aU}. The proof is complete.
\end{proof} 

In the limit $\lambda \to 0$, equation \eqref{Hh} turns into Laplace's, whose
solutions are harmonic functions; moreover, the assumption about $r$ becomes
superfluous in this case. Thus, letting $\lambda \to 0$ in Theorem~1 leads to an
improved formulation of Kuran's theorem because only positive harmonic functions are
involved; see also \cite{AG}.

The reason why $D$ is supposed to be bounded in Theorem~1 is as follows. It is easy
to construct an unbounded domain of finite volume such that $U$ is not integrable
over it, and so boundedness of $D$ allows us to avoid formulating rather complicated
restrictions on the domain.

In the case of sufficiently smooth $\partial D$, the integral $\int_D u (y) \, \D
y$ can be replaced by the flux integral $\int_{\partial D} \partial u / \partial n_y
\, \D S_y$ in the formulation of Theorem~1; here $n$ is the exterior unit normal.
Indeed, we have
\[ \int_D u (y) \, \D y = \lambda^{-2} \int_D \nabla^2 u \, (y) \, \D y = \lambda^{-2}
\int_{\partial D} \partial u / \partial n_y \, \D S_y \, .
\]
This suggests that the following mean flux equality
\begin{equation*}
\frac{\lambda^2 r \, \Gamma (m/2)}{2 \, \Gamma \left( \frac{m}{2} + 1 \right)} \, a
(\lambda r) \, u (x_0) = \frac{1}{|\partial D|} \int_{\partial D} \frac{\partial
v}{\partial n_y} \, \D S_y
\end{equation*}
(cf. formula (31) in \cite{Ku1}) may also characterize the ball of radius $r$
centred at $x_0 \in D$ provided $D$ has a smooth boundary, the point $x_0$ exists
and the equality holds for every solution of equation \eqref{Hh} in $D_r$.

In conclusion we notice that the equality (see \cite{Ku1}, Theorem~8)
\[ m I_{m/2} (\lambda r) \int_{\partial B_r (x)} u (y) \, \D S_y = \lambda r 
I_{(m-2)/2} (\lambda r) \int_{B_r (x)} u (y) \, \D y
\]
holds for every point $x$ belonging to a domain $D \subset \RR^m$ and all $r$ such
that $\overline{B_r (x)} \subset D$ if and only if $u$ is a solution of equation
\eqref{Hh} in $D$. This is analogous to the equality of the mean values over spheres
and balls for harmonic functions. In view of Kosmodem'yanskii's theorem, one might
expect that this equality with $B_r (x)$ changed to $D$ characterizes balls in
$\RR^m$ provided it is valid for every solution of equation \eqref{Hh} in $D$.

\end{document}